\documentclass[12pt,draft]{amsart}
\usepackage[all]{xy}
\usepackage{amsfonts}
\usepackage{amssymb}

\newtheorem{teo}{Theorem}[section]
\newtheorem{lem}[teo]{Lemma}

\newtheorem{cor}[teo]{Corollary}

\theoremstyle{definition}
\newtheorem{dfn}[teo]{Definition}
\newtheorem{rk}[teo]{Remark}
\newtheorem{ex}[teo]{Example}
\def\<{\langle}
\def\>{\rangle}
\def\ss{\subset}

\def\e{\varepsilon}

\def\r{\rho}

\def\t{\tau}

\def\f{{\varphi}}

\def\F{{\Phi}}

\def\C{{\mathbb C}}

\def\Ker{\mathop{\rm Ker}\nolimits}

\def\Tr{\operatorname{Tr}}

\newcommand{\ov}[1]{\overline{#1}}
\newcommand{\til}[1]{\widetilde{#1}}
\newcommand{\wh}[1]{\widehat{#1}}

\newcommand{\Mat}[4]{\left( \begin{array}{cc}
                            #1 & #2 \\
                            #3 & #4
                      \end{array} \right)}

\def\Fix{\operatorname{Fix}}
\def\Rep{\operatorname{Rep}}

\def\cR{{\mathcal R}}
\def\N{{\mathbb N}}
\newcommand\hk[1]{\stackrel\frown{#1}}
\def\cd{\operatorname{cd}}

\begin{document}

\title[ Twisted Burnside theorem for countable groups]
{A  twisted Burnside theorem for countable groups
and Reidemeister numbers}

\author{Alexander Fel'shtyn}
\address{Instytut Matematyki, Uniwersytet Szczecinski,
ul. Wielkopolska 15, 70-451 Szczecin, Poland and Department of Mathematics,
Boise State University,
1910 University Drive, Boise, Idaho, 83725-155, USA}
\email{felshtyn@diamond.boisestate.edu, felshtyn@mpim-bonn.mpg.de}

\author{Evgenij Troitsky}
\thanks{The second author is partially supported by
RFFI Grant  02-01-00574, Grant for the support of
leading scientific schools $\rm HI\! I\! I$-619.203.1
and Grant ``Universities of Russia.''}
\address{Dept. of Mech. and Math., Moscow State University,
119992 GSP-2  Moscow, Russia}
\email{troitsky@mech.math.msu.su}
\urladdr{
http://mech.math.msu.su/\~{}troitsky}

\begin{abstract}
The purpose of the present paper is to prove for finitely
generated groups of type I the following conjecture of
A.~Fel'shtyn and R.~Hill \cite{FelHill}, which is a generalization
of the classical Burnside theorem.

Let $G$ be a countable discrete group, $\phi$ one of its  automorphisms,
$R(\phi)$ the number of $\phi$-conjugacy classes, and
$S(\phi)=\# \Fix (\wh\phi)$ the
number of $\phi$-invariant equivalence classes of irreducible
unitary representations. If one of $R(\phi)$ and
$S(\phi)$ is finite, then it is equal to the other.

This conjecture plays a  important role in the theory of
twisted conjugacy classes (see \cite{Jiang}, \cite{FelshB})
and has very important  consequences
in Dynamics, while its proof needs rather sophisticated  results from
Functional and Non-commutative Harmonic Analysis.

We begin a  discussion of
the general case (which needs another definition of the dual object).
 It will be
the subject of a forthcoming paper.

Some applications
and examples are presented.
\end{abstract}

\maketitle

\tableofcontents

\section{Introduction and formulation of results}

\begin{dfn}
Let $G$ be a countable discrete group and $\phi: G\rightarrow G$ an
endomorphism.
Two elements $x,x'\in G$ are said to be
 $\phi$-{\em conjugate} or {\em twisted conjugate}
iff there exists $g \in G$ with
$$
x'=g  x   \phi(g^{-1}).
$$
We shall write $\{x\}_\phi$ for the $\phi$-{\em conjugacy} or
{\em twisted conjugacy} class
 of the element $x\in G$.
The number of $\phi$-conjugacy classes is called the {\em Reidemeister number}
of an  endomorphism $\phi$ and is  denoted by $R(\phi)$.
If $\phi$ is the identity map then the $\phi$-conjugacy classes are the usual
conjugacy classes in the group $G$.
\end{dfn}

If $G$ is a finite group, then the classical Burnside theorem (see e.g.
\cite[p.~140]{Kirillov})
says that the number of
classes of irreducible representations is equal to the number of conjugacy
classes of elements of $G$.  Let $\wh G$ be the {\em unitary dual} of $G$,
i.e. the set of equivalence classes of unitary irreducible
representations of $G$.

\begin{rk}
If $\phi: G\to G$ is an epimorphism, it induces a map $\wh\phi:\wh G\to\wh G$,
$\wh\phi (\r)=\r\circ\phi$
(because a representation is irreducible if and only if the  scalar operators in the space of
representation are the only ones which commute with all operators of the
representation). This is not the case for a general endomorphism $\phi$,
because $\r \phi$ can be reducible for an irreducible representation $\r$,
and  $\wh\phi$ can be defined only as a multi-valued map.
But
nevertheless we can define the set of fixed points $\Fix \wh\phi$
of $\wh\phi$ on $\wh G$.
\end{rk}

Therefore, by the Burnside's theorem, if $\phi$ is the identity automorphism
of any finite group $G$, then we have
 $R(\phi)=\#\Fix(\wh\phi)$.

To formulate our theorem for the case
of a general endomorphism we first need
an appropriate definition of the $\Fix(\wh\phi)$.

\begin{dfn}
Let $\Rep(G)$ be the space of equivalence classes of
finite dimensional unitary representations of $G$.
Then the corresponding map $\wh\phi_R:\Rep(G)\to \Rep(G)$
is defined in the same way as above: $\wh\phi_R(\r)=\r\circ\phi$.

Let us denote by $\Fix(\wh\phi)$ the set of points $\r\in\wh G\ss
\Rep(G)$ such that $\wh\phi_R (\r)=\r$.
\end{dfn}

\begin{teo}[Main Theorem]\label{teo:mainth1}
Let $G$ be a finitely generated discrete group of type {\rm I},
$\phi$ one of its endomorphism,
$R(\phi)$ the number of $\phi$-conjugacy classes, and
$S(\phi)=\# \Fix (\wh\phi)$ the
number of $\phi$-invariant equivalence classes of irreducible
unitary representations. If one of $R(\phi)$ and
$S(\phi)$ is finite, then it is equal to the other.
\end{teo}

Let $\mu(d)$, $d\in\N$, be the {\em M\"obius function},
i.e.
$$
\mu(d) =
\left\{
\begin{array}{ll}
1 & {\rm if}\ d=1,  \\
(-1)^k & {\rm if}\ d\ {\rm is\ a\ product\ of}\ k\ {\rm distinct\ primes,}\\
0 & {\rm if}\ d\ {\rm is\ not\ square-free.}
\end{array}
\right.
$$

\begin{teo}[Congruences for the Reidemeister numbers]\label{teo:mainth3}
Let $\phi:G$ $\to G$ be an endomorphism of a countable discrete group $G$
such that
all numbers $R(\phi^n)$ are finite and let $H$ be a subgroup
 of $G$ with the properties
$$
  \phi(H) \subset H
$$
$$
  \forall x\in G \; \exists n\in \N \hbox{ such that } \phi^n(x)\in H.
$$
If the pair  $(H,\phi^n)$ satisfies the conditions of
Theorem~{\rm~\ref{teo:mainth1}}
for any $n\in\N$,
then one has for all $n$,
 $$
 \sum_{d\mid n} \mu(d)\cdot R(\phi^{n/d}) \equiv 0 \mod n.
 $$
\end{teo}

These theorems were proved previously in a   special case of
Abelian finitely generated plus finite group \cite{FelHill,FelBanach}.

The interest in twisted conjugacy relations has its origins, in particular,
in the Nielsen-Reidemeister fixed point theory (see, e.g. \cite{Jiang,FelshB}),
in Selberg theory (see, eg. \cite{Shokra,Arthur}),
and  Algebraic Geometry (see, e.g. \cite{Groth}).

Concerning some topological applications of our main results, they are already
obtained in the present paper (Theorem~\ref{teo:febo58}).
The congruences give some necessary conditions for the realization problem
for Reidemeister numbers in topological dynamics.
The relations with Selberg theory will be presented in a forthcoming paper.

Let us remark that it is known that the Reidemeister number of an
endomorphism of a finitely generated Abelian group is
finite iff $1$ is not in the
spectrum of the restriction of this endomorphism to the free part of the group
(see, e.g. \cite{Jiang}).  The Reidemeister number
is infinite for any automorphism of a  non-elementary
Gromov hyperbolic group
\cite{FelPOMI}.

\medskip
To make the presentation more transparent we start from
a new approach (E.T.) for Abelian (Section \ref{sec:abelcase}) and compact
(Section \ref{sec:compcase}) groups.
Only after that we develop this approach and
prove the main theorem for finitely generated groups
of type I in Section \ref{sec:typeI}. A discussion of some examples
leading to conjectures is the subject of Section \ref{sec:exampdisc}.
Then we prove the congruences theorem
(Section \ref{sec:congrReiNumbEnd}) and describe
some topological applications (Section \ref{sec:congrReiNumbMap}).

\medskip\noindent
{\bf Acknowledgement.} We would like to thank the Max Planck
Institute for Mathematics in Bonn for its kind support and
hospitality while the most part of this work has been completed.
We are also indebted to MPI and organizers of the
Workshop on Noncommutative Geometry and Number Theory (Bonn,
August 18--22, 2003)
where the results of this paper were presented.

The first author thanks the Research Institute in Mathematical Sciences
in Kyoto for the possibility of the present research during his visit there.

The authors are grateful to
V.~Balantsev,
B.~Bekka,
R.~Hill, V.~Kaimanovich, V.~Ma\-nui\-lov,
A.~Mish\-che\-n\-ko, A.~Rosenberg,
A.~Shtern, L.~Vainerman,
A.~Vershik for helpful discussions.

They are also grateful to the referee for valuable remarks and
suggestions.

\section{Abelian case}\label{sec:abelcase}

Let $\phi$ be an automorphism of an Abelian group $G$.

\begin{lem}
The twisted conjugacy class $H$ of $e$ is a subgroup.
The other ones  are cosets $gH$.
\end{lem}

\begin{proof}
The first statement follows from the equalities
$$
h\phi(h^{-1}) g\phi(g^{-1})=gh \phi((gh)^{-1},\quad
(h\phi(h^{-1}))^{-1}=\phi(h) h^{-1}=h^{-1} \phi(h).
$$
For the second statement suppose $a\sim b$, i.e. $b=ha\phi(h^{-1})$. Then
$$
gb=gha\phi(h^{-1})=h(ga)\phi(h^{-1}), \qquad gb\sim ga.
$$
\end{proof}

\begin{lem}
Suppose, $u_1,u_2\in G$, $\chi_H$ is the
characteristic function of $H$ as a set. Then
$$
\chi_H(u_1 u_2^{-1})=\left\{
\begin{array}{ll}
1,& \mbox{ if } u_1,u_2 \mbox{ are in one coset },\\
0, & \mbox{ otherwise }.
\end{array}
\right.
$$
\end{lem}

\begin{proof}
Suppose, $u_1\in g_1 H$, $u_2\in g_2 H$, hence, $u_1=g_1 h_1$, $u_2=g_2 h_2$.
Then
$$
u_1 u_2^{-1}=g_1 h_1 h_2^{-1} g_2^{-1} \in g_1 g_2^{-1} H.
$$
Thus, $\chi_H(u_1 u_2^{-1})=1$ if and only if
$g_1 g_2^{-1}\in  H$ and $u_1$ and $u_2$
are in the same class. Otherwise it is 0.
\end{proof}

The following Lemma is well known.

\begin{lem}
For any subgroup $H$ the function $\chi_H$ is of positive type.
\end{lem}

\begin{proof}
Let us take arbitrary elements $u_1,u_2,\dots, u_n$ of $G$.
Let us reenumerate them
in such a way that some first are in $g_1 H$, the next ones are  in $g_2 H$,
and so on, till
$g_m H$, where $g_j H$ are different cosets. By the previous Lemma the matrix
$\|p_{it}\|:=\|\chi_H(u_i u_t^{-1})\| $ is block-diagonal with
square blocks formed by
units. These blocks, and consequently  the whole  matrix are  positively
semi-defined.
\end{proof}

\begin{lem}
In the Abelian case characteristic functions of twisted conjugacy classes
belong to
the Fourier-Stieltjes algebra $B(G)=(C^*(G))^*$.
\end{lem}

\begin{proof}
In this case the characteristic functions of twisted conjugacy
classes are the shifts
of the characteristic function of the class $H$ of $e$. Indeed, we have the
following sequence of equivalent properties:
$$
a\sim b,\qquad b=h a\phi(h^{-1})\mbox{ for some } h,\qquad
gb=gh a\phi(h^{-1})\mbox{ for some } h,
$$
$$
gb=h ga\phi(h^{-1})\mbox{ for some } h,\qquad
ga\sim gb.
$$
Hence, by Corollary (2.19) of \cite{Eym}, these characteristic functions are in
$B(G)$.
\end{proof}

Let us remark that there exists a natural isomorphism (Fourier transform)
$$
u\mapsto\wh u,\qquad C^*(G)=C_r^*(G)\cong C(\wh G),\qquad
\wh g(\rho):=\rho(g),
$$
(this is a  number because irreducible representations of an Abelian group
are 1-dimen\-sional). In fact, it is better to look (for what follows) at an  algebra $C(\wh G)$ as an algebra of continuous sections  of a bundle of 1-dimensional matrix algebras.
over $\wh G$.

Our characteristic functions, being in $B(G)=(C^*(G))^*$ in this case,
are mapped to the functionals on $C(\wh G)$  which, by the
Riesz-Markov-Kakutani theorem,
are measures on $\wh G$. Which of these measures  are invariant under
the induced (twisted) action of $G$ ? Let us remark,
that an invariant non-trivial
functional gives rise to at least one invariant space -- its kernel.

Let us remark, that convolution under the Fourier transform becomes
point-wise multiplication. More precisely, the twisted action, for example,
is defined as
$$
g[f](\rho)=\rho(g)f(\rho)\rho(\phi(g^{-1})),\qquad \rho\in\wh G,\quad
g\in G,\quad f\in C(\wh G).
$$

There are 2 possibilities  for the twisted action of $G$ on the representation
 algebra $A_\rho \cong \C$ : 1) the linear span of the orbit of $1 \in A_\rho$
is equal to all $A_\rho$, 2) and  the opposite case (the action is trivial).

The second case means that the space of interviewing operators between
 $A_\rho$ and $A_{\wh \phi \rho}$ equals $\C$, and  $\rho$ is a fixed point of
the action $\wh \phi:\wh G\to \wh G$. In the first case this is the opposite
situation.

If we have a finite number of such fixed
points, then the space of twisted invariant measures is just the space of
measures
concentrated in these points. Indeed, let us describe the action of $G$
on measures in
more detail.

\begin{lem}
For any  Borel set $E$ one has
$g[\mu](E)=\int_E g[1]\,d\mu$.
\end{lem}

\begin{proof}
The restriction of measure to any Borel set commutes with the action of  $G$,
since the last is point wise on $C(\wh G)$.
For any  Borel set $E$ one has
$$
g[\mu](E)= \int_E 1\, dg[\mu]= \int_E g[1] \,d\mu.
$$
\end{proof}

Hence, if $\mu$ is twisted invariant, then for any  Borel set $E$ and any
$g\in G$ one has
$$
\int_E (1-g[1])\,d\mu =0.
$$

\begin{lem}
Suppose, $f\in C(X)$, where $X$ is a compact Hausdorff space, and $\mu$
is a regular Borel measure on $X$, i.e. a functional on $C(X)$.
Suppose, for any
Borel set $E\subset X$ one has $\int_E f\, d\mu=0$.
Then $\mu (h)=0$ for any $h\in C(X)$
such that $f(x)=0$ implies $h(x)=0$. I.e. $\mu$ is concentrated off the
interior
of $\rm supp\,f$.
\end{lem}

\begin{proof}
Since the functions of the form $fh$ are dense in the space of the
refered to above  $h$'s,
it is sufficient to verify the statement for $fh$. Let us choose an
arbitrary $\e>0$
and a simple function $h'=\sum\limits_{i=1}^n a_i \chi_{E_i}$
such that $|\mu(fh')-\mu(fh)|<\e$.
Then
$$
\mu(fh')=\sum_{i=1}^n \int_{E_i} a_i f \,d\mu =\sum_{i=1}^n a_i
\int_{E_i}f \,d\mu =0.
$$
Since $\e$ is an arbitrary one, we are done.
\end{proof}

Applying this lemma to a twisted invariant measure $\mu$ and $f=1-g[1]$
we obtain that $\mu$ is concentrated at our finite number of fixed
points of $\wh\phi$, because outside of them $f\ne 0$.

If we have an infinite number of fixed points, then the space is
infinite-dimensional
(we have an infinite number of measures
concentrated in finite number of points, each time different)
and Reidemeister number is infinite as well. So, we are done.

\section{Compact case}\label{sec:compcase}

Let $G$ be a compact group, hence $\wh G$ is a discrete space. Then
 $C^*(G)=\oplus M_i$, where $M_i$ are the matrix algebras of
irreducible representations.
The infinite sum is in the following sense:
$$
C^*(G)=\{ f_i\}, i\in \{ 1,2,3,...\}=\hat G, f_i\in M_i,
\| f_i \| \to 0 (i\to \infty).
$$
When $G$ is finite and $\wh G$ is finite this is exactly Peter-Weyl theorem.

A characteristic function of a twisted class
is a functional on $C^*(G)$. For a  finite group it is evident,   for a
general compact group  it is necessary to verify only the measurability of
the twisted class
with the respect to Haar measure, i.e. that twisted class is Borel.
For a compact $G$,
the twisted conjugacy classes being orbits of twisted action are compact and
hence closed.
Then its complement is open, hence Borel, and the class is Borel  too.

Under the identification it passes to a sequence
$\{ \f_i \} $, where $\f_i$ is a functional on $M_i$ (the properties
of convergence can be formulated, but they play no role at the moment).
The conditions of invariance are the following: for each $\r_i \in\wh G$
one has $g[\f_i]=\f_i$, i.e. for any $a\in M_i$ and any $g\in G$ one has
 $\f_i(\r_i(g) a \r_i(\phi(g^{-1})))=\f_i(a)$.

Let us recall  the following well-known fact.

\begin{lem}\label{lem:funcnamat}
Each functional on matrix algebra has form
$a\mapsto \Tr(ab)$ for a fixed matrix $b$.
\end{lem}

\begin{proof} One has $\dim (M(n,\C))'=\dim (M(n,\C))=n\times n$
and looking
at matrices as at operators in $V$, $\dim V=n$, with base $e_i$,
one can remark that
functionals $a\mapsto \<a e_i, e_j\>$, $i,j=1,\dots, n$,
are linearly  independent. Hence,
any functional takes form
$$
a\mapsto \sum_{i,j} b^i_j \<a e_i, e_j\>=
\sum_{i,j} b^i_j a^j _i = \Tr(ba),\qquad
b:=\|b^i_j\|.
$$
\end{proof}

Now we can study invariant ones:
$$
\Tr(b\r_i(g) a \r_i(\phi(g^{-1})))=\Tr(ba),\qquad \forall\,a,g,
$$
$$
\Tr((b-\r_i(\phi(g^{-1})) b \r_i(g))a)=0,\qquad \forall\,a,g,
$$
hence,
$$
b-\r_i(\phi(g^{-1})) b \r_i(g)=0,\qquad \forall\,g.
$$
Since $\r_i$ is irreducible, the dimension of the space of such $b$ is
1 if $\r_i$ is a fixed
point of $\wh\phi$ and 0 in the opposite case. So, we are done.

\begin{rk}
In fact we are only interested in finite discrete case.
Indeed, for a compact $G$,
the twisted conjugacy classes being orbits of twisted action are
compact and hence closed.
If there is a finite number of them, then are open too.
Hence, the situation is more or less
reduced to a discrete group: quotient by the component of unity.
\end{rk}

\section{Extensions and Reidemeister classes}\label{sec:extens}

Consider a group extension respecting homomorphism $\phi$:
 $$
 \xymatrix{
0\ar[r]&
H \ar[r]^i \ar[d]_{\phi'}&  G\ar[r]^p \ar[d]^{\phi} & G/H \ar[d]^{\ov{\phi}}
\ar[r]&0\\
0\ar[r]&H\ar[r]^i & G\ar[r]^p &G/H\ar[r]& 0,}
$$
where $H$ is a normal subgroup of $G$.
The following argument has partial intersection with
\cite{gowon}.

First of all let us notice that the Reidemeister classes of $\phi$ in $G$
are mapped epimorphically on classes of $\ov\phi$ in $G/H$. Indeed,
$$
p(\til g) p(g) \ov\phi(p(\til g^{-1}))= p (\til g g \phi(\til g^{-1}).
$$
Suppose, $R(\phi)<\infty$. Then the previous remark implies
$R(\ov\phi)<\infty$. Consider a class $K=\{h\}_{\t_g\phi'}$, where
$\t_g(h):=g h g^{-1},$ $g\in G$, $h\in H$. The corresponding equivalence
relation is
\begin{equation}\label{eq:klasstaug}
h\sim \til h h g \phi'(\til h^{-1}) g^{-1}.
\end{equation}
Since $H$ is normal,
the automorphism $\t_g:H\to H$ is well defined.
We will denote
by $K$ the image $iK$ as well. By (\ref{eq:klasstaug})
the shift $Kg$ is a subset of $Hg$ is characterized by
\begin{equation}\label{eq:klasstaugsh}
h g\sim \til h (h g) \phi'(\til h^{-1}).
\end{equation}
Hence it is a subset of $\{hg\}_\phi\cap Hg$ and the partition
$Hg=\cup (\{h\}_{\t_g \phi' }) g$ is a subpartition of
$Hg=\cup ( Hg\cap \{hg\}_\phi).$
\begin{lem}\label{lem:ozenfixp}
Suppose, $|G/H|=N<\infty$. Then
$R(\t_g\phi')\le N R(\phi)$. More precisely, the mentioned
subpartition is not more than in $N$ parts.
\end{lem}

\begin{proof}
Consider the following action of $G$ on itself: $x\mapsto gx\phi(g^{-1}).$
Then its orbits are exactly classes $\{x\}_\phi$. Moreover it maps classes
(\ref{eq:klasstaugsh}) onto each other. Indeed,
$$
\til g \til h (h g) \phi'(\til h^{-1}) \phi(\til g^{-1})=
 {\wh h}  \til g (h g)  \phi(\til g^{-1})
\phi'( {\wh h}^{-1})
$$
using normality of the $H$. This map is invertible
($\til g \leftrightarrow \til g^{-1}$), hence bijection.
Moreover, $\til g$ and $\til g \wh h$, for any $\wh h\in H,$
act in the same way. Or in the other words, $H$ is in the stabilizer
of this permutation of classes (\ref{eq:klasstaugsh}). Hence,
the cardinality of any orbit $\le N$.
\end{proof}

Hence, for any finite $G/H$ the number of classes of the form
(\ref{eq:klasstaugsh}) is finite: it is $\le N  R(\phi)$.

\begin{lem}
Suppose,   $H$ satisfies the following
property: for any automorphism of $H$ with finite Reidemeister
number the characteristic functions of Reidemeister classes of $\phi$
are linear combinations
of matrix elements of some finite number of irreducible
finite dimensional  representations of $H$.
Then the characteristic functions of classes
{\rm (\ref{eq:klasstaugsh})}
are linear combinations
of matrix elements of some finite number of irreducible
 finite dimensional  representations of $G$.
\end{lem}

\begin{proof}
Let $\r_1,\r_2,\dots,\r_k$ be the above  irreducible
representations of $H$, $\r$ its direct sum acting on $V$, and
$\pi$ the regular
(finite dimensional) representation of $G/H$.
Let $\r^I_1,\dots,\r^I_k, \r^I$ be the corresponding
induced representations of $G$. Let the characteristic function of
$K$ be represented under the form $\chi_K(h)=\<\r(h)\xi,\eta\>$.
Let $\xi^I\in L^2(G/H,V)$ be defined by the formulas
$\xi^I(\ov e)=\xi \in V$, $\xi^I(\ov g)=0$ if $\ov g\ne \ov e$.
Define similarly $\eta^I$. Then for $h\in iH$ we have
$$
\r^I(h) \xi^I (\ov g)=\r(s(\ov g) h s(\ov g h)^{-1})\xi (\ov g h)=
\r(h s(\ov g) s(\ov g)^{-1})\xi (\ov g )=
\begin{cases}
\r(h)\xi, & \mbox{ if } \ov g = \ov e, \\
0,  & \mbox{ otherwise. }
\end{cases}
$$
Hence,
$\<\r^I(h) \xi^I,\eta^I\>|_{iH}$ is the characteristic function of
$iK$. Let $u,v\in L^2(G/H)$ be such vectors that $\<\pi(\ov g)u,v\>$
is the characteristic function of $\ov e$. Then
$$
\<(\r^I\otimes \pi)(\xi^I\otimes u,\eta^I\otimes v)\>
$$
is the characteristic function of  $iK$. Other characteristic
functions of classes
{\rm (\ref{eq:klasstaugsh})}
are shifts of this one. Hence matrix elements of the representation
$\r^I\otimes \pi$. It is finite dimensional. Hence it can be decomposed
in a finite direct sum of irreducible representations.
\end{proof}

\begin{cor}[of previous two lemmata]\label{lem:osnoskl}
Under the assumptions  of the previous lemma, the
characteristic functions of Reidemeister classes of $\phi$
are linear combinations
of matrix elements of some finite number of irreducible
finite dimensional representations of $G$.
\end{cor}

\section{The case of groups of type I}\label{sec:typeI}

\begin{teo}\label{teo:proptypeI}
Let $G$ be a discrete group of type {\rm I}. Then
\begin{itemize}
\item \cite[3.1.4, 4.1.11]{DixmierEng}
The dual space $\wh G$ is a $T_1$-topological space.

\item \cite{ThomaInvent}
Any irreducible representation of $G$ is finite-dimensional.

\end{itemize}
\end{teo}

\begin{rk}\label{rk:typeIareHausd}
In fact a discrete group  $G$ is of type I if and only if it has a
normal, Abelian subgroup $M$ of finite index. The dimension of any
irreducible representation of $G$ is at most $[G: M]$
\cite{ThomaInvent}.
\end{rk}

Suppose $R=R(\phi)<\infty$, and let $F\ss L^\infty(G)$ be the
$R$-dimensional space of all twisted-invariant functionals on
$L^1(G)$. Let $K\ss L^1(G)$ be the intersection of kernels of
functionals from $F$. Then $K$ is a linear subspace of $L^1(G)$
of codimension $R$. For each $\rho\in\wh G$ let us denote by
$K_\r$ the image $\rho(K)$. This is a subspace of a
(finite-dimensional) full matrix algebra. Let $\cd_\r$ be
its codimension.

Let us introduce the following set
$$
\wh G_F=\{\r\in\wh G \: |\: \cd_\r\ne 0\}.
$$

\begin{lem}
One has $\cd_\r\ne 0$ if and only if  $\r$ is a fixed point
of $\wh\phi$.
In this case $\cd_\r=1$.
\end{lem}

\begin{proof}
Suppose, $\cd_\r\ne 0$ and let us choose a functional $\f_\r$
on the (finite-dimensional full matrix) algebra $\r(L^1(G))$
such that $K_\r\ss \Ker \f_\r$. Then for the corresponding functional
$\f^*_\r=\f_\r\circ\r$ on $L^1(G)$ one has $K\ss \Ker \f^*_\r$.
Hence, $\f^*_\r\in F$ and is twisted-invariant, as well as $\f_\r$.
Then we argue as in the case of compact group (after Lemma
\ref{lem:funcnamat}).

Conversely, if $\r$ is a fixed point of $\wh \phi$,
it gives rise to a (unique up to scaling) non-trivial
twisted-invariant functional $\f_\r$. Let $x=\r(a)$ be any
element in $\r(L^1(G))$
such that $\f_\r(x)\ne 0$. Then $x\not\in K_\r$, because
$\f^*_\r(a)=\f_\r(x)\ne 0$, while $\f^*_\r$ is a twisted-invariant
functional on $L^1(G)$. So, $\cd_\r\ne 0$.

The uniqueness (up to scaling) of the intertwining operator implies
the uniqueness of the corresponding twisted-invariant functional.
Hence, $\cd_\r=1$.
\end{proof}

Hence,
\begin{equation}\label{eq:gfifix}
\wh G_F = \Fix (\wh \phi).
\end{equation}
>From the property $\cd_\r=1$
one obtains for this (unique up to scaling) functional $\f_\r$:
\begin{equation}\label{eq:kerphirho}
\Ker\f_\r=K_\r.
\end{equation}

\begin{lem}\label{lem:RiGF}
$R=\# \wh G_F$, in particular, the set $\wh G_F$ is finite.
\end{lem}

\begin{proof}
First of all we remark that since $G$ is finitely generated
almost Abelian
(cf. Remark \ref{rk:typeIareHausd}) there is a normal Abelian subgroup
$H$ of finite index invariant under all $\phi$. Hence
we can apply Lemma \ref{lem:osnoskl} to $G$, $H$, $\phi$.
So there is  a finite collection of irreducible
representations of $G$ such that
any twisted-invariant functional is a linear combination of
matrix elements of them, i.e. linear combination of functionals on
them. If each of them gives a non-trivial contribution, it has
to be a twisted-invariant functional on the corresponding matrix
algebra. Hence, by the argument above, these representations belong
to $\wh G_F$, and the appropriate functional is unique up to scaling.
Hence, $R\le S$.

Then we use $T_1$-separation property. More precisely,
suppose some points $\r_1,\dots,\r_s$ belong to $\wh G_F$.
Let us choose some twisted-invariant functionals
$\f_i=\f_{\r_i}$ corresponding to these points as it was
described
(i.e. choose some scaling).
Assume that $\|\f_i\|=1$,
$\f_i(x_i)=1$, $x_i\in\r_i(L^1(G))$.
If we can find $a_i\in L^1(G)$ such that
$\f_i(\r_i(a_i))=\f^*_i(a_i)$
is sufficiently large
and $\r_j(a_i)$, $i\ne j$, are sufficiently
small (in fact it is sufficient $\r_j(a_i)$ to be close enough
to $K_j:=K_{\r_j}$),
then $\f^*_j(a_i)$ are small for $i\ne j$,
and $\f^*_i$ are linear independent and hence,
$s<R$.  This would imply $S:=\# \wh G_F \le R$ is finite.
Hence, $R=S$.

So, the problem is reduced to the search of the above  $a_i$.
Let $d=\max\limits_{i=1,\dots,s} \dim \r_i$.
For each $i$ let $c_i:=\|b_i\|$, where $x_i$ is the unitary equivalence
of $\r_i$ and $\wh\phi \r_i$ and $x_i=\r_i(b_i)$.

Let $c:=\max\limits_{i=1,\dots,s} c_i$
and $\e:=\frac 1{2\cdot s^2\cdot d \cdot c}$.

One can find a positive element $a'_i\in L^1(G)$ such that
$\|\r_i(a'_i)\|\ge 1$ and $\|\r_j(a'_i)\|<\e $ for $j\ne i$.
Indeed, since $\r_i$ can be separated from one point, and hence
from the finite number of points:
$\r_j$, $j\ne i$. Hence, one can find
an element $v_i$ such that $\|\r_i(v_i)\|>1$, $\|\r_j(v_i)\|<1$
for $j\ne i$ \cite[Lemma 3.3.3]{DixmierEng}. The same is true for
the positive element $u_i=v_i^*v_i$.
(Due to density we do not distinguish elements of $L^1$ and $C^*$).
Now for a sufficiently large
$n$ the element $a'_i:=(u_i)^n$ has the desired properties.

Let us take $a_i:=a'_i b_i^*$. Then
\begin{multline}\label{eq:ozendiag}
\f^*_i(a_i)=\Tr(x_i \r_i(a_i))=\Tr (x_i \r_i(a'_i) \r_i(b_i)^*)=
\Tr (x_i \r_i(a'_i) x_i^*)=\\
=\Tr(x_i \r_i(a'_i) (x_i)^{-1})
=\Tr( \r_i(a'_i))\ge \frac 1{\dim\r_i}
\ge \frac 1d.
\end{multline}
For $j\ne i$ one has
\begin{equation}\label{eq:ozenvnediag}
\|\f^*_j(a_i)\|=\|\f_j(\r_j(a'_i b_i^*))\|\le c_i \cdot \e.
\end{equation}
Then the $s\times s$ matrix $\F=\f^*_j(a_i)$ can be decomposed into the
sum of the diagonal matrix $\Delta$ and off-diagonal $\Sigma$. By
(\ref{eq:ozendiag}) one has
$\Delta\ge \frac 1d$. By (\ref{eq:ozenvnediag}) one has
$$
\|\Sigma\|\le s^2 \cdot c_i \cdot \e \le s^2 \cdot c \cdot
\frac 1{2\cdot s^2\cdot d \cdot c}= \frac 1{2d}.
$$
Hence, $\F$ is non-degenerate and we are done.
\end{proof}

Lemma \ref{lem:RiGF} together with (\ref{eq:gfifix})
completes the proof of
Theorem \ref{teo:mainth1} for automorphisms.

We need the
following additional observations
for the proof of Theorem \ref{teo:mainth1} for
a general endomorphism
(in which (3) is false for
infinite-dimensional representations).

\begin{lem}\label{lem:endomfin}
\begin{enumerate}
\item If $\phi$ is an epimorphism, then $\wh G$ is $\wh\phi_R$-invariant.

\item For any $\phi$ the set $\Rep(G)\setminus \wh G$ is $\wh\phi_R$-invariant.

\item The dimension of the space of intertwining operators between
$\r\in\wh G$ and $\wh\phi_R(\r)$ is equal to $1$ if and only if
$\r\in \Fix(\wh\phi)$. Otherwise it is $0$.

\end{enumerate}
\end{lem}

\begin{proof}
(1) and (2):
This follows from the characterization of irreducible representation
as that one for which the centralizer of $\r(G)$ consists exactly of
scalar operators.

(3) Let us decompose $\wh\phi_R(\r)$ into irreducible ones. Since
$\dim H_\r=\dim H_{\wh\phi(\r)}$ one has only 2 possibilities:
$\r$ does not appear  in $\wh\phi(\r)$ and the intertwining number is $0$,
otherwise $\wh\phi_R(\r)$ is equivalent to $\r$. In this case
$\r\in \Fix(\wh\phi)$.
\end{proof}

The proof of Theorem \ref{teo:mainth1} can be now repeated
for the general endomorphism
with the new definition of $\Fix(\wh\phi)$. The item (3) supplies us
with the necessary property.

\section{Examples and their discussion}\label{sec:exampdisc}

The natural candidate for the dual object to be used  instead of
$\wh G$ in the case when the different notions of the dual do not
coincide (i.e. for groups more general than type I one groups) is
the so-called quasi-dual $\hk G$, i.e. the set of
quasi-equivalence classes of factor-representations (see, e.g.
\cite{DixmierEng}). This is a usual object when we need a sort of
canonical decomposition for regular representation or group
$C^*$-algebra. More precisely, one needs the support $\hk G_p $
of the Plancherel measure.

Unfortunately the following example shows that this is not the case.

\begin{ex} Let $G$ be a non-elementary Gromov hyperbolic group.
As it was shown by Fel'shtyn \cite{FelPOMI} with the
help of geometrical methods,
for any automorphism $\phi$ of $G$ the Reidemeister number $R(\phi)$
 is infinite. In particular this is true for free group in two generators
$F_2$. But the support $(\hk {F_2})_p $ consists of one point
(i.e. regular representation is factorial).
\end{ex}

The next hope was to exclude from this dual object the $II_1$-points
assuming that they always give rise to an infinite number of twisted
invariant functionals.
But this is also wrong:
\begin{ex}\label{ex:FGW} (an idea of Fel'shtyn realized in \cite{gowon})
Let $G=(Z\oplus Z)\rtimes_\theta  Z $ be the semi-direct product by a
hyperbolic action $\theta(1)=\Mat 2111$. Let
$\phi $
be an automorphism of $G$ whose restriction to  $Z$ is $-id$ and restriction
to  $Z\oplus Z$ is $\Mat 01{-1}0$.
Then $R(\phi)=4$, while the space $\hk G_p $
consists of a single $\rm II_1$-point once again (cf. \cite[p.~94]{ConnesNCG}).
\end{ex}

These examples show that powerful methods of the decomposition
theory do not work for more general classes of groups.

On the other hand Example \ref{ex:FGW}
disproves the old conjecture of Fel'shtyn
and Hill \cite{FelHill} who supposed that the Reidemeister numbers of a injective endomorphism
 for groups of exponential growth are always infinite.
More precisely, this group is amenable and of exponential growth.
Together with some
calculations for concrete groups which are too routine
to be included in the present paper,
this allow us to formulate the following question.

\noindent{\bf Question.}
Is the
Reidemeister number $R(\phi)$ infinite
for any automorphism $\phi$ of (non-amenable) finitely generated group
$G$ containing $F_2$ ?
\smallskip

In this relation  the following example seems to be interesting.
\begin{ex}\cite{FelGon}
For amenable and non-amenable Baumslag-Solitar groups Reidemeister
numbers are always infinite.
\end{ex}

For Example \ref{ex:FGW} recently we have found 4 fixed points of
$\wh\phi$ being finite dimensional irreducible
representations. They give rise to
4 linear independent twisted invariant functionals. These
functionals
 can also be
obtained from the regular factorial representation. There  also exist
fixed points (at least one)
that are  infinite dimensional irreducible
representations. The corresponding
functionals are evidently linear dependent with the first 4. This
example will be presented in detail in a forthcoming paper.

\section{Congruences for Reidemeister numbers of endomorphisms}
\label{sec:congrReiNumbEnd}

 \begin{lem}[\cite{Jiang}]\label{lem:felb7}
For any endomorphism $\phi$ of a group $G$
and any $x\in G$ one has $\phi(x)\in \{x\}_\phi$.
 \end{lem}

 \begin{proof}
$\phi(x)=x^{-1} x \phi(x)$.
 \end{proof}

The following lemma is useful for calculating Reidemeister numbers.
\begin{lem}\label{lem:felb33}
Let $\phi:G\to G$ be any endomorphism
 of any group $G$, and let $H$ be a subgroup
 of $G$ with the properties
$$
  \phi(H) \subset H
$$
$$
  \forall x\in G \; \exists n\in \N \hbox{ such that } \phi^n(x)\in H.
$$
Then
$$
 R(\phi) = R(\phi_H),
$$
 where $\phi_H:H\to H$ is the restriction of $\phi$ to $H$.

\end{lem}

\begin{proof}
Let $x\in G$. Then there is $n$ such that $\phi^n(x)\in H$.
By Lemma~\ref{lem:felb7} it is known that $x$ is $\phi$-conjugate
 to $\phi^n(x)$.
This means that the $\phi$-conjugacy class $\{x\}_\phi$
 of $x$ has non-empty intersection with $H$.

Now suppose that $x,y\in H$ are $\phi$-conjugate,
 i.e. there is a $g\in G$ such that
 $$gx=y\phi(g).$$
We shall show that $x$ and $y$ are $\phi_H$-conjugate,
 i.e. we can find a $g\in H$ with the above property.
First let $n$ be large enough that $\phi^n(g)\in H$.
Then applying $\phi^n$ to the above equation we obtain
 $$ \phi^n(g) \phi^n(x) = \phi^n(y) \phi^{n+1}(g). $$
This shows that $\phi^n(x)$ and $\phi^n(y)$ are $\phi_H$-conjugate.
On the other hand, one knows by Lemma \ref{lem:felb7} that $x$ and $\phi^n(x)$ are
 $\phi_H$-conjugate, and $y$ and $\phi^n(y)$ are $\phi_H$ conjugate,
 so $x$ and $y$ must be $\phi_H$-conjugate.

We have shown that the intersection with $H$ of a
 $\phi$-conjugacy class in $G$ is a $\phi_H$-conjugacy class
 in $H$.
Therefore, we have a map
$$
\begin{array}{cccc}
 Rest : & \cR(\phi) & \to & \cR(\phi_H)\\
        & \{x\}_\phi & \mapsto & \{x\}_\phi \cap H
\end{array}
$$
It is evident that it has the two-sided inverse
$$
 \{x\}_{\phi_H} \mapsto \{x\}_\phi.
$$
Therefore $Rest$ is a bijection and $R(\phi)=R(\phi_H)$.
\end{proof}

\begin{cor}
 Let $H=\phi^n(G)$. Then $R(\phi) = R(\phi_H)$.
\end{cor}

Now we pass to the proof of Theorem~\ref{teo:mainth3}.

\begin{proof}
>From Theorems \ref{teo:mainth1} and  Lemma~\ref{lem:felb33}
it follows immediately that for every $n$
$$
 R(\phi^n)=\#\Fix\left[ \wh{\phi_H}^n:\wh{H}\to\wh{H} \right].
$$

Let $P_n$ denote the number of periodic points of $\wh{\phi_H}$
of least period $n$. One obtains immediately
$$
 R(\phi^n)=\#\Fix\left[ \wh{\phi_H}^n \right] = \sum_{d\mid n} P_d.
$$
Applying M\"obius' inversion formula, we have,
$$
 P_n = \sum_{d\mid n} \mu(d) R(\phi^{\textstyle\frac nd}).
$$
On the other hand, we know that $P_n$ is always divisible
by $n$, because $P_n$ is exactly $n$ times the number of
 $\wh{\phi_H}$-orbits in $\wh{H}$ of cardinality $n$.

Now we pass to the proof of
Theorem \ref{teo:mainth3} for general endomorphisms.

>From Theorem \ref{teo:mainth1}, Lemma~\ref{lem:felb33}
it follows immediately that for every $n$
$$
 R(\phi^n)=R(\phi_H^n)=\# \Fix(\wh{(\phi^n_H)}_R|_{\wh H})
$$

Let $P_n$ denote the number of periodic points of
$(\wh{\phi_H})_R|_{\wh H}$
of least period $n$. One obtains by Lemma \ref{lem:endomfin}~(2)
$$
 R(\phi^n)=\# \Fix(\wh{(\phi^n_H)}_R|_{\wh H})
 = \sum_{d\mid n} P_d.
$$
The proof can be completed as in the case of automorphisms.
\end{proof}

\section{Congruences for  Reidemeister  numbers of a continuous map}
\label{sec:congrReiNumbMap}
Now we pass to the formulation of the topological counterpart
of the main  theorems.
  Let $X$  be a connected, compact
 polyhedron and $f:X\rightarrow X$  be a continuous map.
 Let $p:\til {X}\rightarrow X$ be the universal cover of $X$
 and $\til{f}:\til{X}\rightarrow \til{X}$ a lifting
of $f$, i.e. $p\circ\til{f}=f\circ p$.
Two liftings $\til{f}$ and $\til{f}^\prime$ are called
{\em conjugate} if there is a element$\gamma$ in the deck transformation group
 $\Gamma\cong\pi_1(X)$
such that $\til{f}^\prime = \gamma\circ\til{f}\circ\gamma^{-1}$.
The subset $p(\Fix(\til{f}))\subset \Fix(f)$ is called
{\em the fixed point class of $f$ determined by the lifting class $[\til{f}]$}.
Two fixed points $x_0$ and $x_1$ of $f$ belong to the same fixed point class iff
 there is a path $c$ from $x_0$ to $x_1$ such that $c \cong f\circ c $
(homotopy relative to endpoints).
This fact can be considered as an equivalent definition of a non-empty fixed point class.
 Every map $f$  has only finitely many non-empty fixed point classes, each a compact
 subset of $X$.

The number of lifting classes of $f$ (and hence the number
of fixed point classes, empty or not) is called the {\em Reidemeister number} of $f$,
which is denoted by $R(f)$.
This is a positive integer or infinity.

Let a specific lifting $\til{f}:\til{X}\rightarrow\til{X}$ be fixed.
Then every lifting of $f$ can be written in a unique way
as $\gamma\circ \til{f}$, with $\gamma\in\Gamma$.
So  the elements of $\Gamma$ serve as  ''coordinates'' of
liftings with respect to the fixed $\til{f}$.
Now,  for every $\gamma\in\Gamma$,  the composition $\til{f}\circ\gamma$
is a lifting of $f$ too;  so there is a unique $\gamma^\prime\in\Gamma$
such that $\gamma^\prime\circ\til{f}=\til{f}\circ\gamma$.
This correspondence $\gamma \rightarrow  \gamma^\prime$ is determined by
the fixed $\til{f}$, and is obviously a homomorphism.

\begin{dfn}
The endomorphism $\til{f}_*:\Gamma\rightarrow\Gamma$ determined
by the lifting $\til{f}$ of $f$ is defined by
$$
  \til{f}_*(\gamma)\circ\til{f} = \til{f}\circ\gamma.
$$
\end{dfn}

We shall identify $\pi=\pi_1(X,x_0)$ and $\Gamma$ in the following way.
Choose base points $x_0\in X$ and $\til{x}_0\in p^{-1}(x_0)\subset \til{X}$
once and for all.
 Now points of $\til{X}$ are in 1-1 correspondence with homotopy classes of paths
in $X$ which start at $x_0$:
for $\til{x}\in\til{X}$ take any path in $\til{X}$ from $\til{x}_0$ to
 $\til{x}$ and project it onto $X$;
conversely,  for a path $c$ starting at $x_0$, lift it to a path in $\til{X}$
which starts at $\til{x}_0$, and then take its endpoint.
In this way, we identify a point of $\til{X}$ with
a path class $\langle c \rangle $ in $X$ starting from $x_0$. Under this identification,
$\til{x}_0= \langle e \rangle $ is the unit element in $\pi_1(X,x_0)$.
The action of the loop class $\alpha =  \langle a \rangle \in\pi_1(X,x_0)$ on $\til{X}$
is then given by
$$
\alpha =  \langle a \rangle   :  \langle c \rangle \rightarrow \alpha \cdot  c
= \langle a\cdot c\rangle.
$$
Now we have the following relationship between $\til{f}_*:\pi\rightarrow\pi$
and
$$
f_*  :  \pi_1(X,x_0) \longrightarrow \pi_1(X,f(x_0)).
$$

\begin{lem}
Suppose $\til{f}(\til{x}_0) = \langle w \rangle$.
Then the following diagram commutes:
$$
\xymatrix{
  \pi_1(X,x_0) \ar[r]^{f_*} \ar[dr]_{\til f_*}&   \pi_1(X,f(x_0)) \ar[d]^{w_*} \\
                &       \pi_1(X,x_0)}
$$
where $w_*$ is the isomorphism induced by the path $w$.
\end{lem}
In other words, for every  $\alpha = \langle a \rangle \in\pi_1(X,x_0)$, we have
$$
\til{f}_*(\langle a\rangle )=\langle w(f \circ a)w^{-1}\rangle.
$$
\begin{rk}
In particular, if
 $ x_0 \in p(\Fix( \til f))$ and $ \til x_0 \in \Fix( \til f)$, then $ \til{f}_*=f_*$.
\end{rk}

\begin{lem}[see, e.g. \cite{Jiang}]\label{lem:sovpalgtop}
Lifting classes of $f$ (and hence  fixed point classes, empty or not)
 are in 1-1 correspondence with
 $\til{f}_*$-conjugacy classes in $\pi$,
the lifting class $[\gamma\circ\til{f}]$ corresponding
to the $\til{f}_*$-conjugacy class of $\gamma$.
We therefore have $R(f) = R(\til{f}_*)$.
\end{lem}

We shall say that the fixed point class $p(\Fix(\gamma\circ\til{f}))$,
which is labeled with the lifting class $[\gamma\circ\til{f}]$,
{\it corresponds} to the $\til{f}_*$-conjugacy class of $\gamma$.
Thus $\til{f}_*$-conjugacy classes in $\pi$ serve as ''coordinates'' for fixed
point classes of $f$, once a fixed lifting $\til{f}$ is chosen.

Using Lemma \ref{lem:sovpalgtop} we may apply the Theorem \ref{teo:mainth3}
to the Reidemeister numbers of continuous maps.

\begin{teo}\label{teo:febo58}
Let $f:X\to X$ be a continuous map of a compact polyhedron $X$
such that all numbers $R(f^n)$ are finite.
Let $f_*:\pi_1(X)\to \pi_1(X)$ be an induced endomorphism of the group
 $\pi_1(X)$  and let $H$ be a subgroup of $\pi_1(X)$ with the properties
\begin{enumerate}
\item   $f_*(H) \subset H$,

\item $ \forall x\in \pi_1(X) \; \exists n\in \N \hbox{ such that } f_*^n(x)\in H$.
\end{enumerate}
If the couple $(H,f_*^n)$ satisfies the conditions of
Theorem~{\rm~\ref{teo:mainth1}}
for any $n\in\N$,
then one has for all $n$,
 $$
 \sum_{d\mid n} \mu(d)\cdot R(f^{n/d}) \equiv 0 \mod n.
 $$
\end{teo}


\providecommand{\bysame}{\leavevmode\hbox to3em{\hrulefill}\thinspace}
\providecommand{\MR}{\relax\ifhmode\unskip\space\fi MR }
\providecommand{\MRhref}[2]{%
  \href{http://www.ams.org/mathscinet-getitem?mr=#1}{#2}
}
\providecommand{\href}[2]{#2}

\end{document}